\input amssym.def
\input amssym
\magnification=1200
\parindent0pt
\hsize=16 true cm
\baselineskip=13  pt plus .2pt
$ $

\def\Z{\Bbb Z}

\def\F {{\cal F}}
\def\G {{\cal G}}

\centerline {\bf Lifting finite groups of outer automorphisms of free
groups, surface groups}

\centerline {\bf and their abelianizations}

\bigskip \bigskip

\centerline {Bruno P. Zimmermann}

\medskip

\centerline {Universit\`a degli Studi di Trieste} \centerline {Dipartimento di
Matematica e Informatica} \centerline {34100 Trieste, Italy} \centerline
{zimmer@units.it}

\bigskip \bigskip

Abstract. {\sl We discuss the possibility of lifting finite subgroups, and in
particular finite cyclic subgroups, with respect to the canonical projections
between automorphism and outer automorphism groups of free groups, surface
groups and their abelianizations.}

\bigskip \bigskip

{\bf 1. Introduction}

\medskip

For a group $G$, denote by $Aut \, G$ its automorphism group and by
$Out \, G = Aut \, G/Inn \, G$ its outer automorphism group (automorphisms
modulo inner automorphisms). For a group homomorphism
$\alpha:G \to H$, we say that a subgroup $U$ of $H$ {\it lifts} to $G$ if
there is an injection $\iota:U \to G$ such that $\alpha \circ \iota = id_U$.

\medskip

Let $F_n$ denote the free group of rank $n$, and $\pi_1\F_g$ the fundamental
group of a closed orientable surface $\F_g$ of genus $g$. We consider the
natural projections

$$\alpha:Aut \, F_n \to Out \; F_n,$$
$$ \beta:Out \, F_n \to GL(n,\Z),$$
$$\gamma:Aut \, \pi_1\F_g \to Out \; \pi_1\F_g,$$
$$\delta:Out_+\pi_1\F_g \to Sp(2g,\Z)$$

where $\beta$ and $\gamma$ are obtained by abelianization of $F_n$ and
$\pi_1\F_g$; by $Sp(2g,\Z) \subset SL(2g,\Z)$ we denote the symplectic group
(see [6]), and by $Out_+\pi_1\F_g$ the subgroup of index two of $Out \,
\pi_1\F_g$ induced by orientation-preserving diffeomorphisms of the surface
$\F_g$. It is well-known that the kernels of these four surjections are
torsionfree: this is clear for $\alpha$ and $\gamma$, for $\beta$ and $\delta$
see e.g. [18].

\medskip

The main result is the following

\bigskip

{\bf Theorem}  {\sl Let $n > 2$, $g \ge 2$. For each of the projections
$\alpha, \beta,\gamma$ and $\delta$, there exist finite cyclic subgroups of the
target groups which do not lift.}

\bigskip

We note that $Out \; F_2 \cong GL(2,\Z)$, and $Out_+\pi_1\F_1 \cong SL(2,\Z) =
Sp(2,\Z)$.

\bigskip

{\bf Corollary}  {\sl For $n > 2$ and $g \ge 2$, the projections $\alpha,
\beta,\gamma$ and $\delta$ do not have right inverses.}

\bigskip

Compare also [5,Remark 3 in the introduction] and [8,Theorem 2]. We note that
all target groups $Out \; F_n$, $GL(n,\Z)$, $Out \; \pi_1\F_g$ and $Sp(2g,\Z)$
have torisonfree subgroups of finite index (they are virtually torsionfree),
and there is the more subtle question if such a torsionfree subgroup of finite
index does lift.

\medskip

In the following sections, we consider separately the four cases, commenting
also on the possibility of lifting other types of finite subgroups, in
particular those of maximal order.  Some of the proofs use classical
results; the most technical case is that of $\beta$ which we consider last.

\medskip

We note that, by the positive solution of the Nielsen realization problem,
every finite subgroup $G$ of $Out \; \pi_1\F_g$ can be realized by an action of
$G$ as a group of homeomorphisms of the surface $\F_g$; similarly, every finite
subgroup $G$ of $Out \; F_n$ can be realized by an action of $G$ on a finite
graph with fundamental group $F_n$ (see [20], or the survey [18]). See [3] for
a classification of the finite subgroups of $Out_+\pi_1\F_2$ and
$Out_+\pi_1\F_3$, and [19] for the finite subgroups of $Out \, F_3$.

\bigskip

{\bf 2. The case of $\delta:Out_+\pi_1\F_g  \to Sp(2g,\Z)$}

\medskip

By the positive solution of the Nielsen realization problem ([10]), any finite
subgroup of $Out_+\pi_1\F_g$ can be realized by a finite subgroup of
diffeomorphisms of the surface $\F_g$, and then also, choosing an appropriate
structure of $\F_g$ as a Riemann surface, by a finite group of automorphisms of
the Riemann surface. We note that the Nielsen realization problem for finite
cyclic and solvable groups is a classical result, see [15] for the history of
the problem (see also [17] for the solvable case).

\medskip

Let $a_1,b_1, \ldots,a_g,b_g$ denote a standard symplectic basis of the first
homology  $\Z^{2g} = (\pi_1\F_g)_{ab}$ of the surface $\F_g$ (see [6,chapter
V.3]). Choose any nontrivial symplectic automorphism of finite order of the
subgroup generated by $a_1$  and $b_1$ (the possible orders are 2,3,4 and 6),
and extend it to a symplectic automorphism of the same order of $\Z^{2g}$ by
the identity on the remaining generators $a_2,b_2, \ldots,a_g,b_g$. By
[6,Theorem V.3.3], for $g>2$ this symplectic automorphism is not induced by a
periodic automorphism of a Riemann surface of genus $g$, and hence it does not
lift to $Out_+\pi_1\F_g$.

\medskip

Concerning the case $g=2$, note that $Sp(4,\Z)$ has a subgroup $SL(2,\Z) \times
SL(2,\Z)$, and hence a cyclic subgroup $\Z_3 \times \Z_4$ of order 12; such a
subgroup does not lift to $Out_+\pi_1\F_2$ since, by a result of Wiman, the
maximal order of an orientation-preserving diffeomorphism of $\F_g$ is $4g+2$
(see [16, 4.14.27]).

\bigskip

{\bf Comments}. The maximal possible order of finite subgroups of
$Out_+\pi_1\F_g$ is well-known:

\bigskip

{\bf Theorem 1} ([17])  {\sl For $g>1$, the order of any finite subgroup of
$Out_+\pi_1\F_g$ is bounded above by $84(g-1)$.}

\medskip

The first proof of this has been given in [17,Satz 5.3], in an equivalent
algebraic formulation, as a consequence of a generalized Riemann-Hurwitz
formula; of course, Theorem 1 follows also from the subsequent solution of the
Nielsen realization problem [10] and the classical Riemann-Hurwitz formula.

\medskip

We do not know the maximal order of finite subgroups of the symplectic group
$Sp(2g,\Z)$. However, $Sp(2g,\Z)$ has a finite subgroup $U$ of order $12^gg!$
(permutations of the pairs of generators $a_i,b_i$ and dihedral groups of
order 12 for each pair), and by Theorem 1 these subgroups do not lift.  By
the result of Wiman mentioned above, the maximal order of an
orientation-preserving diffeomorphism of $\F_g$ is $4g+2$; since for almost
all values of $g$ the group  $U$ has finite cyclic subgroups of larger
orders, this gives many cyclic subgroups of $Sp(2g,\Z)$ which do not lift to
$Out_+\pi_1\F_g$.

\bigskip

{\bf 3. The case of $\gamma:Aut \, \pi_1\F_g \to Out \; \pi_1\F_g.$}

\medskip

Since the center of $\pi_1\F_g$ is trivial, every finite subgroup $U$ of
$Out \; \pi_1\F_g$ determines an  extension, unique up to equivalence,
$$1 \to \pi_1\F_g \to E \to U \to 1$$

which is effective (no element of
$E$ acts by conjugation trivially on $\pi_1\F_g$); note that the extension
splits if and only if $U$ lifts to $Aut \, \pi_1\F_g$. Conversely, any such
effective extension defines a subgroup $U$ of $Out
\; \pi_1\F_g$.

\medskip

By lifting to the universal covering, every {\it free} cyclic group
$\Z_n$ of diffeomorphisms of a surface $\F_g$ of genus $g>1$ defines a
torsionfree extension
$$1 \to \pi_1\F_g \to E \to \Z_n \to 1,$$

and hence an inclusion of $\Z_n$ into $Out \; \pi_1\F_g$ which does not lift
to
$Aut \; \pi_1\F_g$. In fact, the group $E$ acts on the universal covering
$\Bbb R^2$ of $\F_g$ (an extension of the universal covering group
$\pi_1\F_g$); if it has torsion, by Smith fixed point theory some element of
prime power order must have fixed points. Examples are the covering
involutions of the orientable 2-fold coverings of the nonorientable surfaces
of genus $g>2$.

\medskip

Alternatively, start with a torsionfree co-compact group $E$ of isometries of
the hyperbolic plane (a surface group) and consider a normal subgroup
$\pi_1\F_g$ with cyclic factor group $\Z_n$; then $\Z_n$ acts freely by
isometries on the surface $\F_g$ defined by the normal subgroup, and the
induced group  $\Z_n$ of outer automorphisms of $\pi_1\F_g$ does not lift to
$Aut \, \pi_1\F_g$.

\bigskip

{\bf Comments}.  Each group $Out \; \pi_1\F_g$ has many noncyclic and
nonabelian finite subgroups (see [3] for $g = 2$ and 3). On the other hand, the
finite subgroups of $Aut \, \pi_1\F_g$ are very special, in fact one has

\bigskip

{\bf Theorem 2} {\sl For $g > 1$, the finite subgroups of $Aut \,
\pi_1\F_g$ are either cyclic or dihedral.}

\medskip

{\it Proof}. A finite subgroup $U$ of $Aut \, \pi_1\F_g$ defines a split
extension
$$1 \to \pi_1\F_g \to E \to U \to 1,$$

so $U$ injects into $E$. By results of Nielsen (see [17,section 2]), the
extension $E$ acts as a group of homeomorphisms of the boundary $S^1$ of the
unit disk (the sphere at infinity of the hyperbolic plane), and it is easy to
see that finite groups of homeomorphisms of $S^1$ are cyclic or dihedral (see
also [17,Lemma 2.1]).

\medskip

Alternatively, one may apply again the solution of the Nielsen realization
problem. By this solution, $U$ can be realized by a group of diffeomorphisms
of the surface $\F_g$, and then also by a group of isometries of a suitable
hyperbolic surface $\F_g$. Lifting to the universal covering of $\F_g$ (the
hyperbolic plane $\Bbb H^2$), this realizes $E$ as a group of isometries of
$\Bbb H^2$, and every finite group of isometries $\Bbb H^2$ is cyclic or
dihedral.

\bigskip

{\bf 4. The case of $\alpha: Aut \, F_n \to Out \; F_n$}

\medskip

A finite subgroup $U$ of $Out \; F_n$ determines an effective extension
$$1 \to F_n \to E \to U \to 1,$$

and the extension splits if and only if $U$ lifts to $Aut \; F_n$.
Conversely, any such extension defines a finite subgroup $U$ of
$Out \; F_n$.

\medskip

Considering extensions
$$1 \to F_n \to F_{n'} \to \Z_m \to 1$$

where also $E$ is a free group $F_{n'}$ (so $(1-n)=m(1-n')$), for each $n$ one
easily constructs finite cyclic subgroups $\Z_m$ of $Out \; F_n$ which do not
lift to $Aut \; F_n$ (e.g. for $n'=2$).

\bigskip

{\bf Comments}.  The finite subgroups of maximal order of $Out \; F_n$ and
$Aut
\; F_n$ are given by the following

\bigskip

{\bf Theorem 3} [14] {\sl For $n > 2$, the maximal order of a finite subgroup
of $Out \; F_n$ and $Aut \; F_n$ is $2^nn!$. For $n >3$, up to conjugation
there is a unique subgroup of maximal order, generated by permutations and
inversions of a system of free generators.}

\medskip

The finite subgroups of  $Out \; F_3$ are determined in [19]. For the possible
orders of finite cyclic subgroups of $Out \; F_n$ and $Aut \; F_n$, see [1] or
[13]. The maximal order of finite abelian subgroups of $Out \; F_n$ and $Aut \;
F_n$ is determined in [2] and equal to $2^n$, for $n>3$.

\vfill  \eject

{\bf 5. The case of $\beta: Out \, F_n \to GL(n,\Z)$}

\medskip

Denoting by $e_1,\ldots,e_n$ the standard basis of $\Z^n$, we define an
automorphism $\phi$ of order six of $\Z^n$ by
$$\phi(e_1) = -e_2, \; \phi(e_2) = e_1 + e_2, \; \phi(e_i) = e_i \;\; {\rm
for} \;\; i \ge 3.$$

We will show that, for $n>2$, the cyclic subgroup of $GL(n,\Z)$ generated by
$\phi$ does not lift to $Out \, F_n$.

\medskip

Consider an extension
$$1 \to \Z^n \to \bar E \to \Z_6 \to 1$$

where a generator of $\Z_n$ induces the automorphism
$\phi$ of $\Z^n$.  If $\bar E$ is the semidirect product $\Z^n \ltimes
\Z_6$ then the abelianization of  $\bar E$ is
$\Z_6 \times \Z^{n-2}$, if the extension does not split the abelianization is
$\Z^{n-2}$, $\Z_2 \times \Z^{n-2}$ or $\Z_3 \times \Z^{n-2}$.

\medskip

Suppose that $\phi$ can be lifted to an outer automorphism of $F_n$ of order
six, represented by an automorphism $\psi$ of $F_n$. Then $\psi$ defines an
extension, unique up to equivalence,
$$1 \to F_n \to E \to \Z_6 \to 1,$$

and the abelianization of $E$ is
$\Z_m \times \Z^{n-2}$, $m=1,2,3$ of 6.

\medskip

The group $E$ is a finite effective extension of the free group $F_n$. By [9],
the finite extension $E$ of the free group $F_n$ is isomorphic to the
fundamental group $\pi_1(\Gamma,\G)$ of a finite graph of finite groups
$(\Gamma, \G)$ (the iterated free product with amalgamation and HNN-extension
over the vertex groups, amalgamated over the edge groups of a maximal tree, the
HNN-generators corresponding to the edges in the complement of the chosen
maximal tree). The {\it Euler characteristic} of $E \cong \pi_1(\Gamma,\G)$ or
of the graph of groups $(\Gamma, \G)$ is
$$ \chi(E) = \chi(\Gamma,\G) = \sum 1/|G_v| - \sum 1/|G_e|$$

where the sum is taken over all vertex groups $G_v$ resp. all edge groups
$G_e$ of $(\Gamma,\G)$. The Euler characteristic behaves multiplicatively
under finite extensions, in particular in our situation we have $\chi(F_n) =
1-n = 6\,\chi(E)$, or
$$-\chi(E) = (n-1)/6.$$

Since the kernel $F_n$ of the surjection of $E$ onto $\Z_6$ is torsionfree,
the vertex and edge groups of $E = \pi_1(\Gamma,\G)$ inject into $\Z_6$ and
hence are cyclic groups of orders 1,2,3 or 6; in the following, we shall
call an edge or vertex with associated group
$\Z_m$ an {\it $m$-edge} or an {\it $m$-vertex}.

\medskip

We shall assume that the graph of groups $(\Gamma,\G)$ is {\it reduced},
i.e. has no non-closed edges such that the edge group coincides with one of
the two vertex groups (such an edge can be contracted obtaining a graph of
groups with fewer edges). Denote by $T$ a maximal tree of the underlying
graph $\Gamma$; then $\Gamma -T$ has exactly $n-2$ edges (considering the
abelianization of $E$). Note that any 6-vertex of
$(\Gamma,\G)$ contributes a direct summand
$\Z_6$ or $\Z_2$ to the abelianization of $E$, and any 2-vertex contributes
 a summand $\Z_2$. Also, all 6-edges, 3-edges and 2-edges of
$(\Gamma,\G)$ are closed, and hence $T$ consists only of 1-edges.

\medskip

Suppose that $\Gamma$ has more than one vertex. Then the contribution of $T$
to $-\chi(E)$ is $\ge 0$, and
$-\chi(E) \ge (n-2)/6$ (considering only the contribution of the edges in
$\Gamma - T$). Since
$-\chi(E) = (n-1)/6$ it follows that $\Gamma -T$ has only 6-edges except
maybe for a single 3-edge and, estimating
$-\chi(E)$ from below, one easily obtains a contradiction.

\medskip

Hence $\Gamma$ has exactly one vertex which has to be a 6-vertex. Then
$E$ is a split extension of $F_n$ and $\Z_6$, and the abelianization of $E$
is $\Z_6 \times \Z^{n-2}$. It follows now easily that
$(\Gamma,\G)$ has no 1-edge, and either one 2-edge and $n-3\;$ 6-edges, or
two 3-edges and $n-4\;$ 6-edges.  In both cases, since the unique vertex
group $\Z_6$ survives in the abelianization, there is a nontrivial subgroup
of the vertex group which is central in $E$.  But then the extension $E$ was
not effective which is a contradiction.

\bigskip

{\bf Comments}.  a) The maximal order finite subgroups of $Out \; F_n$ are
given by Theorem 3. The situation for finite subgroups of $GL(n,\Z)$ is more
complicated. It is shown in [7] that, for values on $n$ larger than some
constant, the maximum value of finite subgroups of $GL(n,\Z)$ is again $2^nn!$,
and that the maximal groups are generated by permuations and inversions of the
standard generators of $\Z^n$.  However, for $n = 2,4,6,7,8,9$ and 10 there are
subgroups of larger orders, the Weyl groups of the exceptional Lie groups of
types $G_2$, $F_4$, $E_6$, $E_7$ and $E_8$ (of orders 12, 1152, 51840, 2903040
and 696729600). On the basis of a result of Weisfeiler, Feit gave a complete
classification of the maximal order finite subgroups of $GL(n,\Z)$ (see
[7],[12]; this uses  the classification of the finite simple groups).

\medskip

For the maximal orders of finite cyclic subgroups of $GL(n,\Z)$, see [12] or
[13]. The maximal orders of finite abelian subgroups of $GL(n,\Z)$ are
determined in [7] and are larger than those for $Out \; F_n$ (see [2]).

\medskip

b)  It is more difficult to construct cyclic subgroups of prime order $p$ of
$GL(n,\Z)$ which do not lift to $Out \; F_n$. Any integral representation of
$\Z_p$ can be written as a direct sum of indecomposable representations which
(in the language of [4,section 1]) are either trivial, regular, cyclotomic or
"exotic" (corresponding to a non-principal ideal in a cyclotomic representation
$\Z[\lambda]$ where  $\lambda$ is a primitive $p$ th root of unity and a
generator of $\Z_p$ acts by multiplication with $\lambda$; so in this case, the
representation has a nontrivial ideal class invariant in the ideal class group
of $\Z[\lambda]$). Now any subgroup $\Z_p$ of $Out \; F_n$ can be induced by
the action of $\Z_p$ on a finite graph with fundamental group $F_n$ ([20]), and
it follows from an argument due to Swan (see [4,section 1]), or from
[11,Theorem 15.5] that the induced representation of $\Z_p$ on the
abelianization of the fundamental group is standard (has no exotic
indecomposable summand). On the other hand, if an integer representation of
$\Z_p$ is standard than it is easy to construct an action of $\Z_p$ on a finite
graph (with one global fixed point) which induces this representations (for
each regular summand one takes a bouquet of $p$ circles permuted cyclically by
a generator of $\Z_p$, for each cyclotomic summand a graph with two vertices
and $p+1$ connecting edges permuted cyclically). Hence the following holds

\bigskip

{\bf Theorem 4} {\sl  A cyclic subgroup $\Z_p$ of prime order $p$ of $GL(n,\Z)$
lifts to $Out \; F_n$ if and only if the corresponding integer representation
of $\Z_p$ is standard.}

\medskip

Exotic integer representations of $\Z_p$ do not exist for $p < 23$. On the
other hand, the situation for general cyclic subgroups $\Z_m$ of $GL(n,\Z)$
appears to be rather complicated.

\bigskip

{\bf Problem.}  Which cyclic subgroups $\Z_m$ of $GL(n,\Z)$ lift to $Out \;
F_n$?

\bigskip \bigskip

\centerline {\bf References}

\medskip

\item {[1]} Z.Bao, {\sl Maximum order of periodic outer automorphisms of a
free group.} J. Algebra 224, 437-453 (2000)

\item {[2]} Z.Bao, {\sl Maximum order of finite abelian subgroups in the
outer automorphism group of a rank $n$ free group.} J. Algebra 236, 355-370
(2001)

\item {[3]} S.A.Broughton, {\sl Classifying finite group actions on
surfaces of low genus.} J. Pure Appl. Algebra 69,  233-270  (1990)

\item {[4]} A.L.Edmonds, {\sl Aspects of group actions on 4-manifolds.}
Top. Appl. 31, 109-124  (1989)

\item {[5]} B.Farb, H.Masur, {\sl Superrigidity and mapping class groups.}
Topology 37, 1169-1176  (1998)

\item {[6]} H.M.Farkas, I.Kra, {\sl Riemann surfaces}. Second edition.
Graduate Texts in Mathematics 71, Springer 1991

\item {[7]} S.Friedland, {\sl The maximal orders of finite subgroups in
$GL_n(\Bbb Q)$.} Proc. Amer. Math. Soc. 125, 3519-3526 (1997)

\item {[8]} N.V.Ivanov, J.D.McCarthy, {\sl On injective homomorphisms
between Teichm\"uller modular groups.} Invent. math. 135, 425-486 (1999)

\item {[9]} A.Karras, A.Pietrowski, D.Solitar, {\sl Finite and
infinite extensions of free groups.} J. Austral. Math. Soc. 16, 458-466 (1972)

\item {[10]} S.P.Kerckhoff, {\sl The Nielsen realization problem.}  Ann.
Math. 117,  235-265  (1983)

\item {[11]} R.S.Kulkarni, {\sl Lattices on trees, automorphisms of graphs,
free groups, and surfaces.}   Unpublished manuscript

\item {[12]} J.Kuzmanovich, A.Pavlichenkov, {\sl Finite groups of matrices
whose entries are integers.}  Monthly 109, 173-186 (2002)

\item {[13]} G.Levitt, J.-L.Nicolas, {\sl On the maximum order of torsion
elements in $GL_n(\Bbb Z)$ and $Aut(F_n)$.} J. Algebra 208, 630-642 (1998)

\item {[14]} S.Wang, B.Zimmermann,  {\it  The maximum order finite groups of
outer automorphisms of free groups.}  Math. Z. 216, 83-87 (1994)

\item {[15]} H.Zieschang, {\sl Finite groups of mapping classes of
surfaces}. Lecture Notes in Mathematics 875, Springer 1981

\item {[16]} H.Zieschang, E.Vogt, H.-D.Coldewey, {\sl Surfaces and planar
discontinuous groups}. Lecture Notes in Mathematics 835, Springer 1980

\item {[17]} B.Zimmermann, {\it  Eine Verallgemeinerung der Formel von
Riemann-Hurwitz.}  Math. Ann. 229, 279-288  (1977)

\item {[18]} B.Zimmermann,  {\it A survey on large finite group actions on
graphs, surfaces and 3-manifolds.} Rend. Circ. Matem. Palermo 52, 47-56
(2003)

\item {[19]} B.Zimmermann,  {\it  Finite groups of outer automorphisms of
free groups.}  Glasgow J. Math. 38, 275-282 (1996)

\item {[20]} B.Zimmermann, {\it \"Uber Hom\"oomorphismen n-dimensionaler
Henkelk\"orper und endliche  Erweiterungen von Schottky-Gruppen.}
Comm. Math. Helv. 56, 474-486 (1981)

\bye